\documentclass[12pt,a4paper,reqno,dvips]{article}%
\usepackage{amsfonts}
\usepackage{bm}
\usepackage[all]{xy}
\usepackage{amsmath}
\usepackage{amssymb}
\usepackage{graphicx}
\SelectTips{cm}{}
\xyoption{dvips}
\newtheorem{theorem}{Theorem}

\newtheorem{definition}[theorem]{Definition}

\newtheorem{proposition}[theorem]{Proposition}

\newenvironment{proof}[1][Proof]{\noindent\textbf{#1.} }{\ \rule{0.5em}{0.5em}}
\title{\bf Iterated Differential Forms II: Riemannian Geometry Revisited}
\author{\sc{A.~M.~Vinogradov}\thanks{{\bf e}-{\it mail}: \texttt{vinograd@unisa.it}} and \sc{L.~Vitagliano}\thanks{{\bf e}-{\it mail}: \texttt{luca\_vitagliano@fastwebnet.it}}\\
\small{DMI, Universit\`a degli Studi di Salerno}\\ \small{and INFN, Gruppo collegato di Salerno,}\\
\small{Via Ponte don Melillo, 84084 Fisciano (SA), Italy}}
\begin{document}
\maketitle
\begin{abstract}
A natural extension of Riemannian geometry to a much wider context
is presented on the basis of the iterated form formalism developed
in \cite{vv06a} and an application to general relativity is given.
\end{abstract}
\newpage
In this note the formalism of iterated differential forms
developed  in the preceding note \cite{vv06a} (see also
\cite{vv06b}) is used to decipher the conceptual meaning of
Riemannian geometry in the context of differential calculus over
(graded) commutative algebras. This is indispensable in order to
extend Riemannian geometry to algebraic geometry and Secondary
Calculus and to clarify some aspects in general relativity and,
more generally, in field theory. In particular, we show that a
Levi-Civita-like connection is naturally associated with a generic
covariant 2-tensor whose symmetric part is nondegenerate and on
this basis present new Einstein-like equations in general
relativity. The super-case is briefly sketched. In the sequel the
notation and definitions of \cite{vv06a} (see also
\cite{v72,v81,v89,n03}) are used.

\section{Differential form calculus and connections\label{SecConn}}

Let $\mathcal{G}=(G,\mu)$ be a grading group, $\Bbbk$ be a field
of zero--characteristic. In what follows a graded
$\Bbbk$--algebra will be always assumed associative, unitary and
graded commutative (see \cite{vv06a}). Accordingly, a graded
morphism of algebras will be always assumed unitary. The necessary
facts concerning differential forms over $\mathcal{G}$--graded
$\Bbbk$--algebras we shall need are summarized below along the
lines of \cite{kv98}.

Let $B$ be a $\mathcal{G}$--graded $\Bbbk$--algebra, $\Lambda (B)=\bigoplus_{k}\Lambda^{k}(B)$ be the
$\mathcal{G}\oplus\mathbb{Z}$--graded $\Bbbk$--algebra of
differential forms over $B$ and $d_{B}:\Lambda
(B)\longrightarrow\Lambda(B)$ the exterior differential. If $P$ is
a graded $B$--module, denote by $\mathrm{D}(B,P)$ the totality of
all $P$--valued graded derivations of $B$.  In a natural way
$\mathrm{D}(B,P)$ is a graded $B$--module.

Let $\varphi:B\longrightarrow C$ be a morphism of
$\mathcal{G}$--graded $\Bbbk$--algebras. It allows to supply any
$C$--module $Q$ with a structure of graded $B$--module which will
be denoted by $Q^{\varphi}$. $\varphi$ induces a morphism of
$\mathcal{G}\oplus\mathbb{Z}$--graded $\Bbbk$--algebras
$\Lambda(\varphi):\Lambda(B)\longrightarrow\Lambda(C)$. Denote by
$\Lambda
^{1}(\varphi):\Lambda^{1}(B)\longrightarrow\Lambda^{1}(C)^{\varphi}$
its first
homogeneous component (with respect to the $%
\mathbb{Z}
$--grading). $\Lambda^{1}(\varphi)$ is a morphism of $B$--modules.
Let $X\in\mathrm{D}(B,\Lambda^{l}(C)^{\varphi})$, $l\geq 0$. The insertion of
$X$ into a differential form $\sigma\in\Lambda(B)$, denoted by
$i_{X}\sigma\in\Lambda(C)$, is defined inductively on the degree
of $\sigma$ by the formula
\[
i_{X}(da\wedge\sigma)=X(a)\wedge(\Lambda(\varphi)(\sigma))-(-1)^{X\cdot
a+l}((d_{C}\circ\varphi)(a))\wedge i_{X}\sigma,\quad a\in B.
\]
The insertion operator
$i_{X}:\Lambda(B)\longrightarrow\Lambda(C)^{\varphi}$ is a
derivation of the algebra $\Lambda(B)$ of bi--degree $(|X|,l-1)$.
The Lie derivative along $X$ is defined as
$\mathcal{L}_{X}\stackrel{\mathrm{def}}{=}\lbrack i_{X},d]=
i_{X}\circ d_{B}+(-1)^{l}d_{C}\circ i_{X}:\Lambda
(B)\longrightarrow\Lambda(C)^{\varphi}$. $\mathcal{L}_{X}$ is a
derivation of the algebra $\Lambda(B)$ of bi--degree $(|X|,l)$.

Let $\Lambda^{1}(B)$ be a projective and finitely generated $B$--module. Then,
for any $X,Y\in$ $\mathrm{D}(B,\Lambda(B))$ there exists a unique
$Z\in\mathrm{D}(B,\Lambda(B))$ such that $[\mathcal{L}_{X},\mathcal{L}%
_{Y}]=\mathcal{L}_{Z}$. Put $Z\stackrel{\mathrm{def}}{=}\lbrack
X,Y]^{\mathrm{fn}}$. The bracket
$[\cdot,\cdot]^{\mathrm{fn}}:\mathrm{D}(B,\Lambda(B))\times\mathrm{D}%
(B,\Lambda(B))\ni(X,Y)\longmapsto\lbrack X,Y]^{\mathrm{fn}}\in\mathrm{D}%
(B,\Lambda(B))$ is called Fr\"{o}licher--Nijenhuis (F--N) bracket and supplies
$\mathrm{D}(B,\Lambda(B))$ with a structure of $\mathcal{G}\oplus\mathbb{Z}%
$--graded $\Bbbk$--Lie--algebra.

Now let $\varphi:B_{1}\longrightarrow B_{2}$ be a morphism of $\mathcal{G}%
$--graded $\Bbbk$--algebras. Consider the functor
\[
\mathrm{D}(B_{1},{}\bullet{}^{\varphi}):P\longmapsto\mathrm{D}(B_{1}%
,P^{\varphi}),
\]
$P$ being a $\mathcal{G}$--graded $B_{2}$--module. It is represented by the
graded $B_{2}$--module $\Lambda^{1}(B_{1})\otimes_{B_{1}}B_{2}^{\varphi}$,
i.e., for any graded $B_2$--module $P$ and $X\in
\mathrm{D}(B_{1},P^{\varphi})$ there exists a unique graded $B_{2}$--module
homomorphism $\psi_{X}:\Lambda^{1}(B_{1})\otimes_{B_{1}}B_{2}^{\varphi
}\longrightarrow P$ such that $X=\psi_{X}\circ d_{\varphi},$ where
$d_{\varphi}:B_{1}\longrightarrow\Lambda^{1}(B_{1})\otimes_{B_{1}}%
B_{2}^{\varphi}$ is a natural derivation defined by
\[
d_{\varphi}b\overset{\mathrm{def}}{=}d_{B_{1}}b\otimes1\in\Lambda^{1}%
(B_{1})\otimes_{B_{1}}B_{2}^{\varphi},\quad b\in B_{1}.
\]
Note that $\Lambda^{1}(\varphi)\circ d_{B_{1}}:B_{1}\longrightarrow\Lambda
^{1}(B_{2})$ is a graded derivation, i.e., $\Lambda^{1}(\varphi)\circ
d_{B_{1}}\in\mathrm{D}(B_{1},\Lambda^{1}(B_{2})^{\varphi})$. Denote
$\psi=\psi_{\Lambda^{1}(\varphi)\circ d_{B_{1}}}:\Lambda^{1}%
(B_{1})\otimes_{B_{1}}B_{2}^{\varphi}\longrightarrow\Lambda^{1}(B_{2})$
and put $\Lambda_{\varphi}^{1}(B_{1},B_{2})=\operatorname{im}\psi
\subset\Lambda^{1}(B_{2})$ (or, simply, $\Lambda^{1}(B_{1},B_{2})$
if the context does not allow a confusion).

\begin{definition}
A \emph{connection in} $\varphi$ is a graded derivation $\nabla:B_{2}%
\longrightarrow\Lambda(B_{1})\otimes_{B_{1}}B_{2}^{\varphi}$ such that,
$\nabla(B_{2})\subset\Lambda^{1}(B_{1})\otimes_{B_{1}}B_{2}^{\varphi}$ and
$\nabla\circ\varphi=d_{\varphi}$.
\end{definition}

Let $\nabla$ be a connection in $\varphi$ and $i_{\nabla}:\Lambda^{1}%
(B_{2})\longrightarrow\Lambda^{1}(B_{1})\otimes_{B_{1}}B_{2}^{\varphi}$ the
unique graded $B_{2}$--module homomorphism such that $\nabla=i_{\nabla}\circ
d_{B_{2}}$. By definition $i_{\nabla}$ is a left inverse of $\psi$. This
proves that $\psi$ is injective and $\Lambda^{1}(B_{1})\otimes_{B_{1}}%
B_{2}^{\varphi}\simeq\Lambda_{\varphi}^{1}(B_{1},B_{2})$.

A connection $\nabla$ in $\varphi$ defines a natural
transformation of functors  $\nabla^{\bullet}:$
$\mathrm{D}(B_{1},{}\bullet
{}^{\varphi})\longrightarrow\mathrm{D}(B_{2},{}\bullet{})$ (see
\cite{kv98}). Namely,
\[
\nabla^{P}:\mathrm{D}(B_{1},P^{\varphi})\ni X\longmapsto\nabla_{X}^{P}%
\stackrel{\mathrm{def}}{=}\psi_{X}\circ\nabla\in\mathrm{D}(B_{2},P),
\]
$P$ being a $\mathcal{G}$--graded $B_{2}$--module. By definition
$\nabla_{X}^{P}\circ\varphi=X$.

There are many situations when a new connection can be constructed
on the basis of other ones. Some of them we need  are described below.

Let
$\nabla:B_{2}\longrightarrow\Lambda(B_{1})\otimes_{B_{1}}B_{2}^{\varphi}$
be a connection in $\varphi:B_{1}\longrightarrow B_{2}$ and
$\kappa :B_{2}\longrightarrow B_{2}$ an automorphism of the
algebra $B_{2}$. Obviously,
$\kappa:B_{2}^{\varphi}\longrightarrow B_{2}^{\kappa\circ\varphi}$ is $B_{1}%
$--linear and, so, the map
\begin{equation}
\mathrm{id}_{\Lambda(B_{1})}\otimes\kappa:\Lambda(B_{1})\otimes_{B_{1}}%
B_{2}^{\varphi}\longrightarrow\Lambda(B_{1})\otimes_{B_{1}}B_{2}^{\kappa
\circ\varphi}. \label{tens}%
\end{equation}
is well-defined.
The composition%
\[
\nabla^{\kappa}\stackrel{\mathrm{def}}{=}(\mathrm{id}_{\Lambda(B_{1})}\otimes\kappa)\circ
\nabla\circ\kappa^{-1}:B_{2}\longrightarrow\Lambda(B_{1})\otimes_{B_{1}}%
^{\kappa\circ\varphi}B_{2}%
\]
is a connection in $\kappa\circ\varphi:B_{1}\longrightarrow B_{2}$.

Let $\varphi$ and $\nabla$ be as above. $\nabla$ extends naturally
to a connection in
$\Lambda(\varphi):\Lambda(B_{1})\longrightarrow\Lambda(B_{2})$ as
follows. Consider the bi--differential algebra
$(\Lambda_{2}(B_{1}),d_{1},d_{2})$ of doubly iterated differential
forms over $B_{1}$ (see \cite{vv06a}).

Note that the map
\[
\nu:B_{2}\ni b\longmapsto1\otimes b\in\Lambda(B_{1})\otimes_{B_{1}}%
B_{2}^{\varphi}%
\]
is a morphism of $\mathcal{G}$--graded algebras and $\nabla\in\mathrm{D}%
(B_{2},(\Lambda(B_{1})\otimes_{B_{1}}B_{2}^{\varphi})^{\nu})$. According to
the above construction $\nabla$ extends to a derivation
\[
\mathcal{L}_{\nabla}:\Lambda(B_{2})\longrightarrow\Lambda(\Lambda
(B_{1})\otimes_{B_{1}}B_{2}^{\varphi}).
\]
There exists a canonical isomorphism of graded algebras,
\[
\phi:\Lambda(\Lambda(B_{1})\otimes_{B_{1}}B_{2}^{\varphi})\;\widetilde
{\longrightarrow}\;\Lambda_{2}(B_{1})\otimes_{\Lambda(B_{1})}\Lambda
(B_{2})^{\Lambda(\varphi)},
\]
given by the formula
\[
(\phi\circ d_{\Lambda(B_{1})\otimes B_{2}})(\omega\otimes b)=d_{2}%
\omega\otimes b+\omega\otimes d_{B_{2}}b,\quad\omega\in\Lambda(B_{1}),\,b\in
B_{2}.
\]
Put
\[
\Lambda\nabla\overset{\mathrm{def}}{=}\phi\circ\mathcal{L}_{\nabla}%
:\Lambda(B_{2})\longrightarrow\Lambda_{2}(B_{1})\otimes_{\Lambda(B_{1}%
)}\Lambda(B_{2})^{\Lambda(\varphi)}.
\]
Then $\Lambda\nabla$ is a connection in $\Lambda(\varphi):\Lambda(B_{1}%
)\longrightarrow\Lambda(B_{2})$.

Finally, let $B_{1},B_{2},B_{3}$ be $\mathcal{G}$--graded algebras,
$\varphi:B_{2}\longrightarrow B_{3}$ and $\varphi^{\prime}:B_{1}%
\longrightarrow B_{2}$ be homomorphisms of graded algebras and $\nabla
:B_{3}\longrightarrow\Lambda(B_{2})\otimes_{B_{2}}B_{3}^{\varphi}$,
$\!\nabla^{\prime}:B_{2}\longrightarrow\Lambda(B_{1})\otimes_{B_{1}}%
B_{2}^{\varphi^{\prime}}$ be connections in $\varphi$,
$\varphi^{\prime}$, respectively. Similarly as above (see
(\ref{tens})) consider the homomorphism
\[
\mathrm{id}_{\Lambda(B_{1})}\otimes\varphi:\Lambda(B_{1})\otimes_{B_{1}}%
B_{2}^{\varphi^{\prime}}\longrightarrow\Lambda(B_{1})\otimes_{B_{1}}%
B_{3}^{\varphi\circ\varphi^{\prime}}.
\]
By the universal property of $\Lambda^{1}(B_{2})\otimes_{B_{2}}%
B_{3}^{\varphi}$ there exists a unique homomorphism of graded $B_{3}$--modules
$\overline{\psi}:\Lambda^{1}(B_{2})\otimes_{B_{2}}B_{3}^{\varphi
}\longrightarrow\Lambda(B_{1})\otimes_{B_{1}}B_{3}^{\varphi\circ
\varphi^{\prime}}$ such that $\overline{\psi}\circ d_{\varphi}=(\mathrm{id}%
_{\Lambda(B_{1})}\otimes\varphi)\circ\nabla^{\prime}$. Put
\[
\square\equiv{}\overline{\psi}\circ\nabla:B_{3}\longrightarrow\Lambda
(B_{1})\otimes_{B_{1}}B_{3}^{\varphi\circ\varphi^{\prime}}.
\]
Then $\square$ is a connection in $\varphi\circ\varphi^{\prime}:B_{1}%
\longrightarrow B_{3}$.
\[
\xymatrix@C=17pt@R=32pt{B_2 \ar[r]^-{\nabla^\prime} \ar[d]_-{d_{\varphi}} & \Lambda(B_1)\otimes B_2^{\varphi^{\prime}} \ar[rr]^-{\mathrm{id}\otimes\varphi}       & &  \Lambda(B_1)\otimes B_3^{\varphi \circ \varphi^{\prime}}    \\
\Lambda^1(B_2)\otimes B_3^{\varphi} \ar[rrru]^-{\overline{\psi}} & & &\\
B_3 \ar[u]^-{\nabla} \ar[rrruu]^-{\square} & & &  }
\]

\begin{definition}
$\square$ is called the \emph{composition} of connections $\nabla$
and $\nabla^{\prime}$.
\end{definition}

Abusing the notation we also denote this connection by
$\nabla\circ\nabla^{\prime}$.

\section{Levi-Civita-like connections over smooth algebras\label{Main}}

To simplify the exposition and for pedagogical reasons we shall limit to the
case when the basic smooth algebra $A$ is the smooth function algebra on a
smooth manifold.

Consider an $n$--dimensional smooth manifold $M$, a local chart $(x^{1}%
,\ldots,x^{n})$ on it and put $A\equiv C^{\infty}(M)$. Let $\Lambda
=\Lambda_{1}$ be the $\mathbb{Z}$--graded algebra of (geometric) differential
forms over $M$ and $\Lambda_{k}=\Lambda(\Lambda_{k-1})$ be the $\mathbb{Z}^{k}%
$--graded algebra of $k$--times iterated (geometric) differential forms
over $M$ (see \cite{ks03}, where iterated differential forms were first introduced, and \cite{s06} for an alternative approach to iterated
differential forms on super--manifolds). Denote by
$\kappa:\Lambda_{2}\longrightarrow\Lambda_{2}$ the canonical
involution (see \cite{vv06a,vv06b}).  Recall that $\kappa\circ
d_1\circ\kappa=d_2$, $\kappa\circ d_2\circ\kappa=d_1$ and $d_{1}$ is a natural extension
of the de Rham differential on $M$ to $\Lambda_{2}$. In the
following we shall denote simply by $\omega\omega^{\prime}$
(rather than $\omega\wedge\omega^{\prime}$) the product of the
iterated forms $\omega$ and $\omega^{\prime}$. Since
$\Lambda_{2}^{1}\equiv\Lambda^{1}(\Lambda)$ is a projective and
finitely generated $\Lambda$--module  the F--N bracket in
$\mathrm{D}(\Lambda,\Lambda_{2})$ is well-defined.

The injective $A$--homomorphism
$\iota_{2}:T_{2}^{0}(M)\hookrightarrow \Lambda_{2}$ (see
\cite{vv06a}) allows one to interpret covariant $2$--tensors, in
particular, (pseudo-)metric tensors, as iterated forms. This way
pseudo--Riemannian geometry becomes a subject of the differential
calculus over the algebra of iterated forms $\Lambda_{\infty}$.
This interpretation allows to associate a linear connection in the
tangent bundle of $M$ with a (possibly non--symmetric) $2$--tensor
whose symmetric part is non--degenerate. According to the standard
approach, the Levi-Civita connection associated with a
pseudo-metric is defined implicitly as the unique one satisfying
certain properties. On the contrary, we define it directly by
applying to the considered $2$--tensor a natural operator in the
algebra $\Lambda_{\infty}$. This is one of numerous examples
showing potentialities of the calculus of iterated forms.

\begin{proposition}
\label{Prop1}Let $g\in\Lambda_{2}$ be a second iterated form of bi--degree
$(1,1)$. Then the following two assertions are equivalent:

\begin{enumerate}
\item The map $\mathrm{D}(M)\times\mathrm{D}(M)\ni(X,Y)\longmapsto(i_{Y}%
^{(2)}\circ i_{X}^{(1)})(g)\in A\subset\Lambda$ is $A$--bilinear,
$\kappa(g)=g$ and the $A$--homomorphism $\rfloor_{g}:\mathrm{D}(M)\ni
X\longmapsto i_{X}^{(2)}g\in\Lambda^{1}\subset\Lambda_{2}$ is bijective.

\item $g=\iota_{2}(g^{\prime})$ for a Riemannian or pseudo--Riemmannian metric
$g^{\prime}$ over $M$.
\end{enumerate}
\end{proposition}

Proposition \ref{Prop1} is a corollary of proposition 4 from
\cite{vv06a} and motivates the following new point of view on
what, in reality, the Riemannian geometry is.

\begin{definition}
\label{newR} A \emph{metric} over a smooth algebra $A$, for instance,
$A=C^{\infty}(M)$, is an iterated $1$--form $g\in\Lambda_{2}$ satisfying
hypothesis 1 of proposition \ref{Prop1}.
\end{definition}
Locally a metric $g\in\Lambda_{2}$ looks as $g=g_{\mu\nu}d_{1}x^{\mu}%
d_{2}x^{\nu}$ with $g_{\mu\nu}=g_{\nu\mu}\in A$,\ $1\leq\mu,\nu\leq n$.

This way the Riemannian geometry is embedded in the context of
iterated differential forms and, after that, it is generalized
rather directly to more general $2$--tensors. Namely, let
$\tau\in\Lambda_{2}$ be an iterated form of bi--degree $(1,1)$
such that $\tau=\iota_{2}(\tau^{\prime})$ for a covariant
$2$--tensor
$\tau^{\prime}$ over $M$. Locally $\tau$ looks as $\tau=\tau_{\mu\nu}%
d_{1}x^{\mu}d_{2}x^{\nu}$ with $\tau_{\mu\nu}\in A$,\
$1\leq\mu,\nu\leq n$. $\tau$ can be uniquely written as
$\tau=g+\omega$ where $\kappa(g)=g$ and $\kappa(\omega)=-\omega$.
Clearly, $g=\iota_{2}(g^{\prime})$ and $\omega
=\iota_{2}(\omega^{\prime})$, where $g^{\prime}$ and
$\omega^{\prime}$ are symmetric and skew--symmetric parts of
$\tau^\prime$, respectively. From now on, suppose that $g$ is
nondegenerate, i.e., a pseudo-metric.

Consider the iterated $2$--form $\gamma=-d_{2}d_{1}\tau$. It is of
bi--degree $(2,2)$ and its local expression is
\[
\gamma=\ \partial_{\nu}\partial_{\beta}\tau_{\alpha\mu}d_{1}x^{\beta}%
d_{1}x^{\alpha}d_{2}x^{\mu}d_{2}x^{\nu}+\gamma_{\beta}{}_{\alpha\mu}%
d_{1}x^{\alpha}d_{2}x^{\mu}d_{2}d_{1}x^{\beta}+\tau_{\alpha\beta}d_{2}%
d_{1}x^{\alpha}d_{2}d_{1}x^{\beta}.
\]
with $\partial_{\mu}=\tfrac{\partial}{\partial x^{\mu}}$ and
$\gamma_{\beta}{}_{\alpha\mu}=\partial_{\mu}\tau_{\alpha\beta}%
+\partial_{\alpha}\tau_{\beta\mu}-\partial_{\beta}\tau_{\alpha\mu}$.
Note that if $\omega=0$, then the $\gamma_{\beta}{}_{\alpha\mu}$'s
are doubled first kind Christoffel symbols of $g^{\prime}$. By
this reason $\gamma$ is naturally baptized to be the
\emph{(generalized) Christoffel form}.

Next, the homomorphism $\rceil^{g}=\tfrac{1}{2}\rfloor_{g}^{-1}:\Lambda
^{1}\longrightarrow\mathrm{D}(M)$ extends to $\Lambda$ as an $A$--linear,
$\mathrm{\ D}(A,\Lambda)$--valued derivation of $\Lambda$ of degree $-1$
denoted by $\rceil_{1}^{g}:\Lambda\longrightarrow\mathrm{D}(A,\Lambda)$. It
should be noted that $\mathrm{D}(A,\Lambda)$ is a $\Lambda$--module and
$\rceil_{1}^{g}\in\mathrm{D}(\Lambda,\mathrm{D}(A,\Lambda))$. Recall that with
any $X\in\mathrm{D}(A,\Lambda)$ the insertion operator $i_{X}\in
\mathrm{D}(\Lambda,\Lambda)$ is associated and consider the map
\[
\Delta:\Lambda\ni\sigma\longmapsto\Delta(\sigma)\stackrel{\mathrm{def}}{=} i_{\rceil_{1}%
^{g}(\sigma)}\in\mathrm{D}(\Lambda,\Lambda).
\]
$\Delta$ is an $A$--linear, $\mathrm{D}(\Lambda,\Lambda)$--valued derivation
of $\Lambda$ of degree $-2$. Define the map
\[
\rceil_{2}^{g}:\Lambda_{2}\longrightarrow\mathrm{D}(\Lambda,\Lambda_{2}),
\]
by putting for any $\Omega\in\Lambda_{2}$ and $\omega\in\Lambda$
\[
\rceil_{2}^{g}(\Omega)(\omega)\equiv(-1)^{\Omega\cdot\omega} i_{\Delta(\omega)}^{(2)}\Omega.
\]
Note that $\mathrm{D}(\Lambda,\Lambda_{2})$ possesses a natural $\Lambda_{2}%
$--module structure and $\rceil_{2}^{g}$ is a $\mathrm{D}(\Lambda,\Lambda
_{2})$--valued derivation of $\Lambda_{2}$ of bi--degree $(-2,-1)$. Locally it
looks as
\[
\rceil_{2}^{g}(\Omega)(\omega)=(-1)^{|\Omega|\cdot(0,-1)}\tfrac{1}{2}g^{\mu
\nu}(i_{i_{\partial_{\mu}}}^{(2)}\Omega)(i_{\partial_{\nu}}\omega)\in
\Lambda_{2},
\]
where $\left\Vert g^{\mu\nu}\right\Vert =\left\Vert g_{\mu\nu}\right\Vert
^{-1}$.

\begin{definition}
The $\Lambda_{2}$--valued derivation
$\Gamma=\rceil_{2}^{g}(\gamma)=-
\rceil _{2}^{g}(d_{1}d_{2}\tau)$ of
the algebra $\Lambda$ is called the\emph{ Levi--Civita connection
form} of $\tau$.
\end{definition}

This terminology is motivated by the fact that the local
expression of $\Gamma$ is
\[
\Gamma=(d_{1}d_{2}x^{\alpha}+\Gamma_{\mu}{}_{\beta}^{\alpha}d_{1}x^{\beta
}d_{2}x^{\mu})\,i_{\partial_{\alpha}},
\]
with $\Gamma_{\mu}{}_{\beta}^{\alpha}\equiv\tfrac{1}{2}g^{\alpha\delta}%
\gamma_{\delta}{}_{\beta\mu}$ and in the case when $\omega=0$
$\Gamma_{\rho}{}_{\gamma }^{\alpha}$'s are nothing but Christoffel
symbols of $g^{\prime}$. So, with any covariant $2$--tensor whose
symmetric part is nondegenerate a linear connection in the tangent
bundle is associated (see below).

Now consider the insertion into $2$-iterated forms operator
$i_{\Gamma}^{(2)}\in\mathrm{D}(\Lambda_{2},\Lambda_{2})$
corresponding to the graded form--valued derivation $\Gamma$. The
composition $\kappa\circ i_{\Gamma}^{(2)}\circ\kappa$ is again a
derivation of $\Lambda_{2}$. Put
\[
\kappa(\Gamma)\stackrel{\mathrm{def}}{=}[\kappa\circ i_{\Gamma}^{(2)}\circ\kappa,d_{2}]|_{\Lambda}%
\in\mathrm{D}(\Lambda,\Lambda_{2}).
\]
The graded form--valued derivation $T=\Gamma-\kappa(\Gamma)$ has the following
local expression
\[
T=T_{\mu}{}_{\beta}^{\alpha}d_{1}x^{\beta}d_{2}x^{\mu}\,i_{\partial_{\alpha}%
},\quad T_{\mu}{}_{\beta}^{\alpha}=\Gamma_{\mu}{}_{\beta}^{\alpha}-\Gamma
_{\beta}{}_{\mu}^{\alpha}=3g^{\alpha\delta}\partial_{\lbrack\mu}\omega
_{\beta\delta]}d_{1}x^{\beta}d_{2}x^{\mu}\,i_{\partial_{\alpha}},
\]
where the $\omega_{\beta\delta}$'s are local components of $\omega$, i.e.,
locally $\omega=\omega_{\beta\delta}d_{1}x^{\beta}d_{2}x^{\delta}$, and
indexes in square bracket are skew--symmetrized.

\begin{definition}
The $\Lambda_{2}$--valued derivation $T$ of the algebra $\Lambda$ is called
the torsion of $\tau$.
\end{definition}

The graded F--N square of $\Gamma$%
\[
R=[\Gamma,\Gamma]^{\mathrm{fn}}\in\mathrm{D}(\Lambda,\Lambda_{2})
\]
is well defined and has the following local expression
\[
R=[\Gamma,\Gamma]^{\mathrm{fn}}=R_{\sigma\mu}{}_{\delta}^{\alpha}%
d_{1}x^{\delta}d_{2}x^{\sigma}d_{2}x^{\mu}\,i_{\partial_{\alpha}},
\]
where $R_{\sigma\mu}{}_{\delta}^{\alpha}=\partial_{\sigma}\Gamma_{\mu}%
{}_{\delta}^{\alpha}-\partial_{\mu}\Gamma_{\sigma}{}_{\delta}^{\alpha}%
+\Gamma_{\sigma}{}_{\beta}^{\alpha}\Gamma_{\mu}{}_{\delta}^{\beta}-\Gamma
_{\mu}{}_{\beta}^{\alpha}\Gamma_{\sigma}{}_{\delta}^{\beta}$.

\begin{definition}
The $\Lambda_{2}$--valued derivation $R$ of the algebra $\Lambda$ is called
the Riemann tensor of $\tau$.
\end{definition}

The introduced above concepts of generalized Levi-Civita connection, torsion
and Riemann tensor associated with the $2$--tensor $\tau$ have standard
counterparts. Namely, recall that $\tau=\iota_{2}(\tau^{\prime})$,
$\tau^{\prime}\in T_{2}^{0}(M)$ being a \textquotedblleft
standard\textquotedblright\ covariant $2$--tensor on $M$. The above
computations show that there is a standard connection $\Gamma^{\prime}$ in the
tangent bundle of $M$ associated with $\tau^{\prime}$. Indeed, the covariant
derivative along $X\in\mathrm{D}(M)$ with respect to $\Gamma^{\prime}$
$\nabla_{X}^{\prime}$ is well defined by
\[
(\nabla_{X}^{\prime}Y)(f)=(i_{X}^{(1)}\circ i_{Y}^{(2)}\circ\Gamma)(df)\in
A\subset\Lambda_{2}.
\]
Now, denote by $T^{\prime}$ and $R^{\prime}$ the \textquotedblleft
standard\textquotedblright\ torsion and Riemann curvature tensor of
$\Gamma^{\prime}$, respectively.

\begin{proposition}
For any $X,Y,Z\in\mathrm{D}(M)$ and $f\in A$
\begin{align*}
T^{\prime}(X,Y)(f)  &  =(i_{X}^{(1)}\circ i_{Y}^{(2)}\circ T)(df)\in
A\subset\Lambda_{2}.\\
R^{\prime}(X,Y)(Z)(f)  &  =(i_{Z}^{(1)}\circ i_{X}^{(2)}\circ i_{Y}^{(2)}\circ
R)(df)\in A\subset\Lambda_{2}.
\end{align*}

\end{proposition}

\begin{proof}
Obvious from local expression.
\end{proof}

In new terms geodesics of $\Gamma^{\prime}$ are described as follows. Let
$\tfrac{d}{dt}:C^{\infty}(\mathbb{R})\longrightarrow C^{\infty}(\mathbb{R})$
be the \textquotedblleft time current\textquotedblright\ field. Consider a
smooth curve $\chi:I\longrightarrow M$ on $M$, $I\subset\mathbb{R}$ being an
interval, the algebra $\Lambda_{2}(I)$ of second iterated forms on $I$ and the
pull--back homomorphism $\chi^{\ast}:\Lambda_{2}\rightarrow\Lambda_{2}(I)$.
The \emph{geodesic curvature} vector field (along $\chi$) is the composition
\[
K_{\chi}\overset{\mathrm{def}}{=}i_{d/dt}\circ i_{d/dt}^{(2)}\circ\chi^{\ast
}\circ\Gamma:\Lambda\longrightarrow\Lambda(I).
\]

\begin{proposition}
\label{Prop3} $\chi$ is an affinely parameterized geodesic of the metric
$g^{\prime}$ if and only if $K_{\chi}=0$.
\end{proposition}

This is obvious from the local expression
\[
K_{\chi}(dx^{\mu})=\Gamma_{\rho}{}_{\alpha}^{\mu}\tfrac{d\chi^{\rho}}%
{dt}\tfrac{d\chi^{\alpha}}{dt}+\tfrac{d^{2}\chi^{\mu}}{dt^{2}}.
\]

We now describe covariant derivatives of covariant tensors in new terms. First
of all note that for any $l<k$, $\Lambda(\Lambda_{l})\otimes_{\Lambda_{l}%
}\Lambda_{k}$ is naturally isomorphic to $\Lambda(\Lambda_{l},\Lambda_{k})$.
Therefore, a connection in the inclusion $\Lambda_{l}\subset\Lambda_{k}$ (in
the sense of section \ref{SecConn}) may be understood as a derivation
$\Lambda_{k}\longrightarrow\Lambda(\Lambda_{l},\Lambda_{k})$. In the following
we will adopt this point of view.

Define a new graded form-valued derivation $\nabla=d_{2}-\Gamma:\Lambda
\longrightarrow\Lambda_{2}.\nabla$ has local expression
\[
\nabla=d_{2}x^{\mu}\,(\partial_{\mu}-\Gamma_{\mu}{}_{\beta}^{\alpha}%
\,d_{1}x^{\beta}\,i_{\partial_{\alpha}}).
\]
Note that $\operatorname{im}\nabla\subset\Lambda(A,\Lambda)$.

\begin{proposition}
$\nabla:\Lambda\longrightarrow\Lambda(A,\Lambda)$ is a connection in the
natural inclusion $\varphi:A\hookrightarrow\Lambda$ in the sense of section
\ref{SecConn}.
\end{proposition}

As shown in section \ref{SecConn}, in its turn $\nabla$ induces a connection
$\Lambda\nabla$ in $\Lambda(\varphi):\Lambda\hookrightarrow\Lambda_{2}$ . Note
that $\Lambda(\varphi)=\kappa|_{\Lambda}$. Therefore,
\[
\nabla_{2}\equiv\Lambda(\kappa)\circ\nabla^{\prime}\circ\kappa:\Lambda
_{2}\longrightarrow\Lambda(\Lambda,\Lambda_{2})
\]
is a connection in the natural inclusion $\Lambda\subset\Lambda_{2}$.
Iterating the procedure we get, at the $k$--th step, a connection $\nabla_{k}$
in the inclusion $\Lambda_{k-1}\subset\Lambda_{k}$. Note that the existence of
the connections $\nabla_{k}$ enriches noteworthy the standard calculus in
Riemannian geometry (consequences of this remark will be explored elsewhere).
Denote by $\overset{k}{\nabla}$ the composition
\[
\overset{k}{\nabla}{}=\nabla_{k}\circ\nabla_{k-1}\circ\cdots\circ\nabla
_{2}\circ\nabla_{1}:\Lambda_{k}\longrightarrow\Lambda(A,\Lambda_{k}).
\]
$\overset{k}{\nabla}$ is a connection in the inclusion $A\subset\Lambda_{k}$.
Accordingly, we can lift $A$--linearly a derivation $X\in\mathrm{D}%
(A,\Lambda_{k})$ to a derivation $\overset{k}{\nabla}_{X}=i_{X}^{(k+1)}\circ
{}\overset{k}{\nabla}{}\in\mathrm{D}(\Lambda_{k},\Lambda_{k})$.

\begin{definition}\label{conn}
$\overset{k}{\nabla}$ is called the $k$\emph{--th covariant derivative
associated to }$\Gamma$.
\end{definition}

As an example, we report local expression of $\overset{2}{\nabla}$:
\[
\overset{2}{\nabla}{}=d_{3}x^{\mu}(\partial_{\mu}-\Gamma_{\mu}{}_{\beta
}^{\alpha}\,d_{1}x^{\beta}\,i_{\partial_{\alpha}}-\Gamma_{\mu}{}_{\beta
}^{\alpha}\,d_{2}x^{\beta}\,i_{\partial_{\alpha}}^{(2)}-d_{1}(\Gamma_{\mu}%
{}_{\beta}^{\alpha}\,d_{2}x^{\beta})\,i_{i_{\partial_{\alpha}}}^{(2)}).
\]
The terminology of definition \ref{conn} is motivated by the fact that the standard
covariant derivatives with respect to $\Gamma^{\prime}$ of a
covariant $k$--tensor $s^{\prime}\in T_{k}^{0}(M)$ are naturally
recovered from $\overset{k}{\nabla}$. Namely, we have

\begin{proposition}
Let $s=\iota_{k}(s^{\prime})\in\Lambda_{k}$ and $X,Y_{1},\ldots,Y_{k}%
\in\mathrm{D}(M)$. Then
\[
(\nabla_{X}^{\prime}\,s^{\prime})(Y_{1},\ldots,Y_{k})=(i_{Y_{1}}^{(1)}\circ
{}\cdots\circ i_{Y_{k}}^{(k)}\circ i_{X}^{(k+1)}\circ{}\overset{k}{\nabla}%
{})(s)\in A\subset\Lambda_{k}.
\]

\end{proposition}

It is worth mentioning that, on the contrary to the the classical
Levi-Civita connections, i.e., associated with metrics, the
covariant derivative of a $2$--tensor $\tau$ with respect to the
associated connection is, generally, different from zero, i.e.,
$\overset{2}{\nabla}\tau\neq0$. Namely, let $\Gamma^{g}:\Lambda
_{1}\longrightarrow\Lambda_{2}$ be the (generalized) Levi-Civita
connection of $g$ and $\overset{k}{\nabla}{}^{g}$ the associated
$k$--th covariant derivative. Then,
\[
\overset{2}{\nabla}\tau={}\overset{2}{\nabla}\omega={}\overset{2}{\nabla}%
{}^{g}\omega+T(\omega),
\]
where $T(\omega)$ is the iterated form locally given by
\[
T(\omega)=T_{\mu\lbrack\nu}^{\beta}\omega_{\alpha]\beta}d_{1}x^{\alpha}%
d_{2}x^{\nu}d_{3}x^{\mu}.
\]

\section{Natural equations in general relativity}

The possibility to associate in a canonical way a Levi-Civita-like
connection with a $2$--tensor $\tau$ with non--degenerate
symmetric part suggests an immediate application to general
relativity. Indeed, it is natural to interpret the skew-symmetric
part $\omega$ of $\tau$ as a matter field. Then, it is natural as
well to assume that the space, i.e., the symmetric part $g$ of
$\tau$, is shaped by the matter in such a way that the Ricci
tensor $\mathrm{Ric}(\tau)$ of the associated with $\tau$
connection vanishes. This way one gets the \emph{natural}
equations in general relativity. We use here the word \textquotedblleft natural\textquotedblright{} in order to
stress that in the proposed approach the credit is definitively
given to the naturalness of the background mathematical language,
but  not to the \textquotedblleft physical intuition\textquotedblright{} which could be easily
misleading in this context.

In coordinates these natural equations look as follows. Let $R=R_{\mu\sigma}%
{}_{\delta}^{\alpha}dx^{\delta}\otimes dx^{\sigma}\otimes dx^{\mu}%
\otimes\tfrac{\partial}{\partial x^{\alpha}}$ be the associated
with $\tau$ standard Riemann tensor. The Ricci tensor
$\mathrm{Ric}(\tau)$ of $\tau$ is defined as the $2$--tensor
locally given by
$\mathrm{Ric}(\tau)=R_{\mu\delta}\,dx^{\mu}\otimes dx^{\delta}$,
where $R_{\mu\delta}=R_{\mu\alpha}{}_{\delta}^{\alpha}$,
$\mu,\delta=1,\ldots,n$ (the interpretation of the Ricci tensor
in terms of iterated differential forms will be discussed
separately). Then the natural equations for $\tau$ read
\bigskip$\mathrm{Ric}(\tau)=0$.

Now denote by $\overset{g}{\nabla}$ and $\mathrm{Ric}(g)$ the
covariant differential and the Ricci tensor of the pseudo-metric
$g$, respectively, and put
$\mathrm{Ric}(g)={}\overset{g}{R}_{\mu\delta}\,dx^{\mu}\otimes
dx^{\delta}$. Then in terms of $g$ and $\omega$ the natural
equations for $\tau$ read
\begin{equation}
{}\left.
\begin{array}
[c]{c}%
\overset{g}{R}_{\mu\nu}+\tfrac{9}{16}g^{\gamma\rho}g^{\delta\sigma}%
\partial_{\lbrack\delta}\omega_{\mu\rho]}\partial_{\lbrack\gamma}\omega
_{\nu\sigma]}=0\\
\overset{g}{\nabla}{}^{\gamma}(\partial_{\lbrack\mu}\omega_{\nu\gamma]})=0
\end{array}
\right.  .\label{Einst}%
\end{equation}
System (\ref{Einst}) may be understood as Einstein-like equations
for $g$ in presence of a matter field (represented by $\omega$)
satisfying a suitable constitutive equation. Indeed, it can be
shown that (\ref{Einst}) are Einstein equations in presence of a
perfect irrotational fluid with a definite equation of state
(details will be discussed in a separate publication).

It is worth stressing that the matter in the above natural
equations is treated in a simplified manner. More rich and exact
versions of natural equations in the frames of the developed
approach will be discussed separately.

\section{Levi-Civita connections over smooth super-algebras}

The previous results  are straightforwardly generalized to
supermanifolds. In this section, we shall sketch how Riemannian
supergeometry (see, for instance, \cite{g06} and references
therein) can be developed in the framework of iterated
differential forms along the lines of section \ref{Main}.

Consider a super--manifold $S$, a local chart $(\theta^{1},\ldots,\theta^{n})$
on it and put $A=C^{\infty}(S)$. $A$ is a $%
\mathbb{Z}
_{2}$--graded commutative algebra. Put $\bar{\alpha}=|\theta^{\alpha}|{}\in%
\mathbb{Z}
_{2}$. Let $\Lambda=\Lambda_{1}$ be the $%
\mathbb{Z}
_{2}\mathbb{\oplus Z{}}$--graded algebra of (geometric) differential forms
over $S$ and $\Lambda_{2}=\Lambda(\Lambda)$ be the $\mathbb{{}}%
\mathbb{Z}
_{2}\mathbb{\oplus Z}^{2}$--graded algebra of doubly iterated
geometric differential forms over $S$ (see also \cite{ks03}).

Let $g$ be an even doubly iterated form of bi--degree $(1,1)$ over
$A$ such
that the map $\mathrm{D}(S)\times\mathrm{D}(S)\ni(X,Y)\longmapsto(i_{Y}%
^{(2)}\circ i_{X}^{(1)})(g)\in A\subset\Lambda$ is $A$--bilinear,
$\kappa(g)=g$ and the $A$--homomorphism $\rfloor_{g}:\mathrm{D}(M)\ni
X\longmapsto i_{X}^{(2)}g\in\Lambda^{1}\subset\Lambda_{2}$ is bijective.

\begin{definition}
$g$ is called a \emph{supermetric over }$S$.
\end{definition}

Locally, $g=g_{\mu\nu}d_{1}\theta^{\mu}d_{2}\theta^{\nu}$, $g_{\mu\nu
}=(-1)^{\bar{\mu}\bar{\nu}}g_{\nu\mu}\in A$. Moreover, the matrix $\left\Vert
g_{\nu\mu}\right\Vert $ is invertible in $A$.

Given a metric $g$ over $S$ the derivation $\rceil_{2}^{g}:\Lambda
_{2}\longrightarrow\mathrm{D}_{\Lambda}(\Lambda_{2})$ is defined
exactly as above. Its local expression is
\[
\rceil_{2}^{g}(\Omega)(\omega)=(-1)^{|\Omega|\cdot(\bar{\nu},0,-1)}\tfrac
{1}{2}g^{\mu\nu}(i_{i_{\partial_{\mu}}}^{(2)}\Omega)(i_{\partial_{\nu}}%
\omega)\in\Lambda_{2},
\]
where $\Omega\in\Lambda_{2}$, $\omega\in\Lambda$, $(\bar{\nu},0,-1,)\in%
\mathbb{Z}
_{2}\mathbb{\oplus\mathbb{Z}}^{2}\mathbb{{}}$ is the multi--degree of
$i_{\partial_{\nu}}$, $\partial_{\nu}=\tfrac{\partial}{\partial\theta^{\nu}}$
and $g^{\mu\nu}\in A$, $(-1)^{\bar{\gamma}\cdot\bar{\alpha}}g_{\nu\alpha
}g^{\nu\gamma}=\delta_{\alpha}^{\gamma}$.

\begin{definition}
The $\Lambda_{2}$--valued derivation $\Gamma=\rceil_{2}^{g}(\gamma)=-\rceil
_{2}^{g}(d_{1}d_{2}g)$ of the algebra $\Lambda$ is called the
\emph{Levi-Civita connection form} of $g$.
\end{definition}

The local expression of $\Gamma$ is
\[
\Gamma=(d_{1}d_{2}\theta^{\alpha}+\Gamma_{\mu}{}_{\beta}^{\alpha}d_{1}%
\theta^{\beta}d_{2}\theta^{\mu})\,i_{\partial_{\alpha}},
\]
with $\Gamma_{\mu}{}_{\beta}^{\alpha}=\tfrac{1}{2}g^{\gamma\alpha}%
[(-1)^{\bar{\gamma}\cdot\bar{\gamma}+\bar{\mu}(\bar{\mu}+\bar{\gamma}%
)}\partial_{\mu}g_{\beta\gamma}+(-1)^{\bar{\gamma}\cdot\bar{\gamma}+\bar
{\beta}(\bar{\mu}+\bar{\beta}+\bar{\gamma})}\partial_{\beta}g_{\mu\gamma
}-\partial_{\gamma}g_{\beta\mu}]$. Similarly, the Riemann tensor and the
covariant derivative are introduced. For instance,%
\[
R=[\Gamma,\Gamma]^{\mathrm{fn}}=R_{\gamma\beta}{}_{\alpha}^{\nu}d_{1}%
\theta^{\alpha}d_{2}\theta^{\gamma}d_{2}\theta^{\beta}i_{\partial_{\nu}}%
\]
where
\begin{align*}
R_{\gamma\beta}{}_{\alpha}^{\nu}  &  =(-1)^{\bar{\alpha}(\bar{\beta}%
+\bar{\gamma})+\bar{\beta}(\bar{\beta}+\bar{\nu})}\partial_{\beta}%
\Gamma_{\alpha}{}_{\gamma}^{\nu}-(-1)^{\bar{\beta}(\bar{\beta}+\bar{\gamma
}+\bar{\nu})}\partial_{\gamma}\Gamma_{\alpha}{}_{\beta}^{\nu}\\
&  +(-)^{\bar{\mu}\bar{\beta}+\bar{\alpha}\bar{\gamma}}\Gamma_{\mu}{}_{\beta
}^{\nu}\Gamma_{\alpha}{}_{\gamma}^{\mu}-(-1)^{\bar{\beta}\cdot(\bar{\alpha
}+\bar{\gamma})+\bar{\mu}\bar{\gamma}}\Gamma_{\mu}{}_{\gamma}^{\nu}%
\Gamma_{\alpha}{}_{\beta}^{\mu}.
\end{align*}

\end{document}